\documentclass{amsart}
\numberwithin{equation}{section}

\newtheorem{thm}{Theorem}[section]
\newtheorem{lem}[thm]{Lemma}

\newtheorem{defin}[thm]{Definition}

\begin{document}

\begin{center}
\textbf{{\large {\ A Nonlinear Nonlocal Problem for the Caputo Fractional Subdiffusion Equation }}}\\[0pt]
\medskip \textbf{Ravshan Ashurov$^{1}$, Rajapboy Saparboyev$^{1}$ and Navbahor Nuraliyeva$^{1}$}\\[0pt]
\textit{ashurovr@gmail.com\\[0pt]}
\textit{rajapboy1202@gmail.com\\[0pt]}
\textit{n.navbahor2197@gmail.com\\[0pt]}
\medskip \textit{$^{1}$ Institute of Mathematics, Academy of Science of Uzbekistan}

\end{center}

\textbf{Abstract}: In this paper, we study a time-fractional subdiffusion equation with a nonlinear nonlocal initial condition involving the unknown solution at the final time. The considered problem is formulated using the Caputo fractional derivative of order \(0 < \alpha < 1\), along with homogeneous Dirichlet boundary conditions. The nonlocal initial condition is of the form \( u(x,0) = g(x, u(x,T)) \), where \(g\) is a nonlinear function satisfying a Lipschitz condition. The main challenge arises from the implicit dependence on the unknown final state. Using an explicit representation of the solution in terms of the Green function and applying the Banach fixed point theorem, we establish the existence and uniqueness of a regular solution. We also provide uniform estimates for the Green function and analyze the influence of the Lipschitz constant on solvability.

\vskip 0.3cm \noindent {\it AMS 2000 Mathematics Subject
Classifications} :
Primary ?35R11; Secondary 34A12.\\
{\it Key words}:   Fractional subdiffusion equation, Caputo fractional derivative, non-local boundary condition, a priori estimate, Green’s function, successive approximations, existence of solution.

\section{Introduction}
Consider the following initial-boundary value problem for the fractional subdiffusion equation posed in the domain \( (0,1) \times (0,T] \):
\begin{equation}\label{Nuraliyevaprob1}
\partial_{0t}^{\alpha }u(x,t) - u_{xx}(x,t) = f(x,t),
\end{equation}
subject to the nonlocal initial condition
\begin{equation}\label{NuraliyevaCauchy1}
u(x,0) = g(x, u(x,T)), \quad x \in [0,1],
\end{equation}
and the Dirichlet boundary conditions
\begin{equation}\label{Nuraliyevadirac1}
u(0,t) = u(1,t) = 0, \quad t \in [0,T],
\end{equation}
where \( \partial_{0t}^{\alpha} \) denotes the Caputo fractional derivative of order \( 0 < \alpha < 1 \) with respect to the time variable \( t \), and \( f(x,t) \), \( g(x,w) \) are given functions.

\begin{defin}\label{Nuraliyevadefen1}
A function $u(t)\in C([0,1]\times [0,T])$ with the properties $\partial_{0t}^{\alpha }u(x,t),u_{xx}(x,t)\in C((0,1)\times (0,T])$ and satisfying conditions \eqref{Nuraliyevaprob1}--\eqref{Nuraliyevadirac1} is called a regular solution of the non-local problem \eqref{Nuraliyevaprob1}--\eqref{Nuraliyevadirac1}.
\end{defin}

Many scientists have studied various non-local conditional problems for subdiffusion equations.

If \( g = \psi(x) \), then such initial-boundary value problems have been studied in papers such as \cite{NuraliyevaAshMuh} and \cite{NuraliyevaMuh1}.

It should be noted that various problems have been considered even when the function \( g \) is linear, as can be seen in \cite{NuraliyevaAshFay}, \cite{Nuraliyevanav2}, and \cite{Nuraliyevanav1}.

The subdiffusion equation with the condition
\[
u(\xi) = \delta u(0) + \varphi, \quad 0 < \xi \leq T, \quad \delta = \text{const}
\]
(instead of condition~(\ref{NuraliyevaCauchy1})) is studied in detail in \cite{NuraliyevaAshFay}. In that work, the values of the parameter that ensure the existence and uniqueness of the solution are determined. In other cases, the authors establish certain orthogonality conditions for \( f(t) \) and \( \varphi \) that guarantee the existence of a solution; however, uniqueness is not ensured in those cases.

In \cite{Nuraliyevanav2} and \cite{Nuraliyevanav1}, the authors employed linear non-local conditions depending on three parameters to solve equation~(\ref{Nuraliyevaprob1}). 

 The work \cite{Nuraliyevanav2} involves a pointwise non-local condition of the form:
\[
\alpha u(0) + \beta u(T) = \gamma u(\xi) + \varphi,
\]
while the paper \cite{Nuraliyevanav1}  includes an integral non-local condition:
\[
\alpha u(0) + \beta u(T) + \gamma \int_0^T u(\eta)\, d\eta = \varphi.
\]
Both problems are formulated in terms of the Caputo fractional derivative of order \( 0 < \alpha < 1 \), and the elliptic part is represented by a self-adjoint positive operator in a separable Hilbert space. The authors establish existence and uniqueness theorems for the solutions of both problems. Furthermore, they identify sufficient conditions that guarantee the uniqueness of the solution. The influence of the parameters \( \alpha \), \( \beta \), and \( \gamma \) on the existence and uniqueness of the solutions is thoroughly analyzed.

In all of the aforementioned works, the Fourier method was employed. 
In contrast, the present study differs by incorporating a nonlinear non-local condition, 
and an unconventional approach is used to solve the problem.

\section{Preliminaries}

In this section, we introduce the definitions of the fractional integral and derivative, present an auxiliary problem essential for proving the main theorems, and state necessary lemmas.

Let \( h(t) \) be a function defined on the interval \([a,b]\), and let \(\sigma > 0\). The Riemann--Liouville fractional integral of order \(\sigma\) is defined by (see \cite{NuraliyevaKilbas}):
\begin{equation}\label{Nuraliyevadef0}
 D_{0t}^{-\sigma} h(t) = \frac{1}{\Gamma(\sigma)} \int_0^t (t - \tau)^{\sigma - 1} h(\tau) \, d\tau,   
\end{equation}
provided that the right-hand side exists point-wise. As usual, $\Gamma(\sigma)$ is Euler's gamma function.

The Caputo fractional derivative of order $0<\rho <1$  for the function  $h(t)$ is defined as (see, e.g., \cite{NuraliyevaKilbas}, p. 92):
$$
\partial_{0t}^{\alpha} h(t) =\frac{1}{\Gamma
(1-\rho)}\int\limits_0^t\frac{h'(\xi)}{(t-\xi)^{\rho}} d\xi,
\quad t>0,
$$
provided that the right-hand side exists point-wise.

For $0 < \rho \leq 1$ and an arbitrary complex number $\mu$, by $
E_{\rho, \mu}(z)$ we denote the Mittag--Leffler function with two parameters (see, e.g., \cite{NuraliyevaKilbas}, p. 56):
$$
E_{\rho, \mu}(z)= \sum\limits_{n=0}^\infty \frac{z^n}{\Gamma(\rho
n+\mu)}.
$$
If the parameter $\mu =1$, then we have the classical Mittag-Leffler function: $ E_{\rho}(z)= E_{\rho, 1}(z)$. Note also $E_{1, 1}(z)=E_{1}(z)=e^z$.

The Beta function \( B(a, b) \) is defined for real numbers \( a > 0 \), \( b > 0 \) by the improper integral

\begin{equation}\label{Nuraliyevabeta}
 B(a, b) = \int_0^1 t^{a-1}(1 - t)^{b-1} \, dt.  
\end{equation}

The Beta function is related to the Gamma function by the identity

\begin{equation}\label{Nuraliyevabetagamma}
B(a, b) = \frac{\Gamma(a)\Gamma(b)}{\Gamma(a + b)}.
\end{equation}

Consider the following auxiliriary problem to find the function  \( u(x,t) \)  
\begin{equation}\label{Nuraliyevaprob10}
\partial_{0t}^{\alpha} u(x,t) - u_{xx}(x,t) = f(x,t),
\end{equation}
subject to the initial and boundary conditions
\begin{equation}\label{NuraliyevaCauchy10}
u(x,0) = \varphi(x), \quad 0 \leq x \leq 1,
\end{equation}
\begin{equation}\label{NuraliyevaDirac}
u(0,t) = 0, \quad u(1,t) = 0, \quad 0 \leq t \leq T.
\end{equation}
where $\partial_{0t}^{\alpha}$ is the Caputo fractional derivative of order $0 <\alpha<1$ in the time variable, $f(x,t),\varphi(x)$ are given functions.

The solution of the problem \eqref{Nuraliyevaprob10}--\eqref{NuraliyevaDirac} is defined in the same manner as in Definition~\ref{Nuraliyevadefen1}.

A solution to problem \eqref{Nuraliyevaprob10}--\eqref{NuraliyevaDirac} can be represented by the formula 
\begin{equation}\label{NuraliyevaGrin}
u(x,t) = \int_0^1 \varphi(\xi) D_{0t}^{\alpha-1} G(x,t,\xi,0) \, d\xi + \int_0^t \int_0^1 G(x,t,\xi,\tau) f(\xi,\tau) \, d\xi d\tau,
\end{equation}
where \( G(x,t,\xi,\tau) \) is the Green function corresponding to the problem. This  function is constructed in Remark 6.1 of \cite{NuraliyevaMamchuyev} and has the form
\begin{equation}\nonumber
G(x,t,\xi,\tau) = \sum_{m=-\infty}^{+\infty} \left[ P(2m + x - \xi, t - \tau) - P(2m + x + \xi, t - \tau) \right], \quad 0 < \xi < 1, \quad 0 < \tau < t \leq T,
\end{equation}
where
\[
P(x,t) = \frac{t^{\frac{\alpha}{2} - 1}}{2} e^{1, \frac{\alpha}{2}}_{1, \frac{\alpha}{2}} \left(-|x| t^{-\frac{\alpha}{2}} \right),
\]
and
\[
e^{\mu,\delta}_{\gamma,\beta}(z) = \sum_{n=0}^\infty \frac{z^n}{\Gamma(\gamma n + \mu) \Gamma(\delta - \beta n)}, \quad \gamma > 0, \quad \gamma > \beta,
\]
is the Wright-type function (see \cite{NuraliyevaPshu}, p. 23).

\begin{lem}\label{eu}[see,\cite{NuraliyevaPshu}, p. 46] \begin{itemize}
    \item If $\delta\ge 0$, $\beta\in(0,1)$, then $e^{1,\delta}_{1,\beta}(-x)>0$ for any positive $x$.
    \item If $\delta\ge \beta$, then when $x>0$ the function $e^{1,\delta}_{1,\beta}(-x)$ is strictly decreasing. 
\end{itemize}
\end{lem}

\begin{lem}\label{e1}[see \cite{NuraliyevaPshu}, p. 47]
If \(\delta \geq 1\), $\beta\in(0,1)$ then for any positive \(x\), the inequalities
\[
0<e^{1,\delta}_{1,\beta}(-x) \leq \frac{1}{\Gamma(\delta)}e^{-x^{\frac{1}{1-\beta}} \beta^{\frac{\beta}{1-\beta}} (1-\beta)}.
\]
\end{lem}

\begin{lem}[see \cite{NuraliyevaPshu}, p. 49]\label{e10} Let $\delta<1$, $\beta\in(0,1)$, then for any positive $x$ and $t$, $a\in (0,x)$, $\xi\in[\beta,1]$, $\omega\in(\frac{1}{2},\min\{1,\frac{1}{2\beta} \})$ is the inequality true
\[
\bigg|t^{\delta-1}e^{1,\delta}_{1,\beta}(-\frac{x}{t^{\beta}})\bigg|\leq \frac{1}{\beta \pi}\big(C_{\beta}(\xi,\omega)\big)^{\frac{\delta-1}{\beta}}\Gamma(\frac{1-\delta}{\beta})(x-a)^{\frac{\delta-1}{\beta}}e^{1,1}_{1,\beta}(-\frac{a}{t^{\beta}}),
\] 
where 
\[
C_{\beta}(\xi,\omega)=(1-\beta)\bigg( \frac{-cos\omega\pi}{\xi-\beta}\bigg)^{\frac{1-\beta}{1-\beta}}\bigg( \frac{cos\beta\omega\pi}{1-\xi}\bigg)^{\frac{1-\beta}{1-\beta}}.
\]
\end{lem}
Using Lemma \ref{e1}, we obtain the following estimate:
\begin{equation}\label{Ge}
\bigg|t^{\delta-1}e^{1,\delta}_{1,\beta}(-\frac{x}{t^{\beta}})\bigg|\leq \frac{1}{\beta \pi}\big(C_{\beta}(\xi,\omega)\big)^{\frac{\delta-1}{\beta}}\Gamma(\frac{1-\delta}{\beta})(x-a)^{\frac{\delta-1}{\beta}}e^{-x^{\frac{1}{1-\beta}} \beta^{\frac{\beta}{1-\beta}} (1-\beta)}=C(x-a)^{-\gamma_{0}} e^{-kx^{\gamma}},
\end{equation}
where 
\[
C=\frac{1}{\beta \pi}\big(C_{\beta}(\xi,\omega)\big)^{\frac{\delta-1}{\beta}}\Gamma(\frac{1-\delta}{\beta}), \quad \gamma_{0},\gamma, k>0.
\]
To estimate the Green function \( G(x,t,\xi,\tau) \), we decompose it into three parts:
\[
G(x,t,\xi,\tau) = G_{-1}(x,t,\xi,\tau) +G_{0}(x,t,\xi,\tau)+ G_1(x,t,\xi,\tau),
\]
where
\[
G_{-1}(x,t,\xi,\tau) = \sum_{m=-\infty}^{-1} \left[ P(2m + x - \xi, t-\tau) - P(2m + x + \xi, t-\tau) \right],
\]
\[
G_{0}(x,t,\xi,\tau)=P(x-\xi,t-\tau)-P(x+\xi,t-\tau),
\]
\[
G_1(x,t,\xi,\tau) = \sum_{m=1}^\infty \left[ P(2m + x - \xi, t-\tau) - P(2m + x + \xi, t-\tau) \right].
\]

By substituting \( m = -n \) in \( G_{-1}(x,t,\xi,\tau) \), we rewrite
\[
G_{-1}(x,t,\xi,\tau) = \sum_{n=1}^\infty \left[ P(-2n + x - \xi, t-\tau) - P(-2n + x + \xi, t-\tau) \right].
\]

Applying the estimate \eqref{Ge} to \( G_1 (x,t,\xi,\tau)\), we have
\[
|G_{-1}(x,t,\xi,\tau)| \leq \sum_{n=1}^\infty \left| P(-2n + x - \xi, t-\tau) \right| + \left| P(-2n + x + \xi, t-\tau) \right|
\]
\[
\leq\sum_{n=1}^{\infty} C|-2n+x-\xi|^{-\gamma_{0}} e^{-k|x-\xi-2n|^{\gamma}}+C|-2n+x+\xi|^{-\gamma_{0}} e^{-k|x+\xi-2n|^{\gamma}}\leq C.
\]
Similarly, for \( G_1 (x,t,\xi,\tau)\), the same bound holds:
\[
|G_1(x,t,\xi,\tau)|\leq C.
\]
Next, we consider an estimate for the function \( G_0(x, t, \xi, \tau) \).

\[
|G_{0}(x,t,\xi,\tau)|\leq |P(x-\xi,t-\tau)|+|P(x+\xi,t-\tau)|.
\]
According to Lemma~\ref{eu}, we obtain the following estimate:
\[
|P(x-\xi,t-\tau)|+|P(x+\xi,t-\tau)|\leq \frac{(t-\tau)^{\frac{\alpha}{2} - 1}}{2} \bigg|e^{1, \frac{\alpha}{2}}_{1, \frac{\alpha}{2}} \left(-|x-\xi| (t-\tau)^{-\frac{\alpha}{2}} \right)\bigg|
\]
\[
+\frac{(t-\tau)^{\frac{\alpha}{2} - 1}}{2} \bigg|e^{1, \frac{\alpha}{2}}_{1, \frac{\alpha}{2}} \left(-|x+\xi| (t-\tau)^{-\frac{\alpha}{2}} \right)\bigg|\leq C (t-\tau)^{\frac{\alpha}{2}-1}.
\]

Hence, the Green function satisfies
\begin{equation}\label{NuraliyevaGrinasim}
|G(x,t,\xi,\tau)| \leq C |t - \tau|^{\frac{\alpha}{2} - 1}, \quad 0 < \tau < t \leq T.
\end{equation}

\section{Main result}
\begin{thm}\label{Nuraliyevatheo2}
Assume that:
\begin{enumerate}
    \item \( f(x,t)\in C([0,1] \times [0,T]) \);
    \item \( g(\xi,v) \in C([0,1] \times \mathbb{R}) \) satisfies the Lipschitz condition
    \[
    |g(\xi, v_1) - g(\xi, v_2)| \leq L |v_1 - v_2| \quad \text{for all } \xi \in [0,1], \, v_1,v_2 \in \mathbb{R};
    \]
    \item \( L < \dfrac{\Gamma\left(1 - \frac{\alpha}{2}\right)}{C\Gamma\left(\frac{\alpha}{2}\right) T^{1 - \frac{\alpha}{2}}} \)  for $C$, which is in (\ref{NuraliyevaGrinasim}).
\end{enumerate}
Then the problem \eqref{Nuraliyevaprob1}--\eqref{Nuraliyevadirac1} admits a unique solution.
\end{thm}
From \eqref{NuraliyevaGrin} we have
\begin{equation}\nonumber
u(x,t) = \int_0^1 \varphi(\xi) D^{\alpha-1}_{0t} G(x,t,\xi,0) d\xi + \int_0^t \int_0^1 G(x,t,\xi,\tau) f(\xi,\tau) d\xi d\tau.
\end{equation}

Setting \(u(x,T) = v(x)\), we obtain:
\begin{equation}\label{Nuraliyevavx}
v(x) = \int_0^1 D^{\alpha-1}_{0T} G(x,T,\xi,0) g(\xi, v(\xi)) d\xi + F(x),
\end{equation}
where
\begin{equation}\nonumber
F(x) = \int_0^T \int_0^1 G(x,T,\xi,\tau) f(\xi,\tau) d\xi d\tau.
\end{equation}

We prove the existence of a solution to the equation in the space \( C[0,1] \). This space consists of all continuous functions on the interval \([0,1]\), equipped with the following norm:
\[
\|v\| = \max_{x \in [0,1]} |v(x)|.
\]

Define the operator \( A: C[0,1] \to C[0,1] \) as follows:
\[
(Av)(x):= \int_0^1 D^{\alpha - 1}_{0T} G(x,T,\xi,0) g(\xi, v(\xi)) \, d\xi + F(x).
\]

We denote the right-hand  side of the equation \eqref{Nuraliyevavx} by $Av(x)$, where $A$ is the corresponding nolinear operator. To apply the Banach fixed point theorem we have to prove the following:

(1) if $v(x)\in C[0, 1]$, then $Av(x)\in C[0, 1]$;

(2) for any $v_1(x),v_2(x)\in C[0,1]$ one has
\[
\max\limits_{t\in [0,1]}||Av_1(x)-Av_2(x)||\leq \delta \max\limits_{t\in [0,t_1]} ||v_1(x)-v_2(x)||, \,\, \delta<1.
\]

We now show that the operator $A$ satisfies condition (1). To this end, we successively make use of the properties of the operator $D^{\alpha-1}$, the upper bound for $|F(x)|$, the Lipschitz condition for the function $g$, and the estimate (\ref{NuraliyevaGrinasim}) for $G(x,T,\xi,0)$:

$$
|v(x)| \le\int_0^1 |D^{\alpha-1}_{0T} G(x,T,\xi,0) g(\xi, v(\xi))| d\xi + |F(x)|
$$
$$
\le\frac{1}{\Gamma(1-\alpha)}\left| \int_0^1 \int_0^T (T-\eta)^{-\alpha} G(x,\eta,\xi,0) g(\xi, v(\xi))d\eta d\xi\right| + \frac{2C}{\alpha}T^{\frac{\alpha}{2}}\max_{[0,1]\times[0,T]}|f(x,t)|
$$
$$
\le\frac{1}{\Gamma(1-\alpha)}\left| \int_0^1 \int_0^T (T-\eta)^{-\alpha} G(x,\eta,\xi,0) g(\xi, v(\xi))d\eta d\xi\right| + \frac{2C}{\alpha}T^{\frac{\alpha}{2}}\max_{[0,1]\times[0,T]}|f(x,t)|
$$
$$
\le\frac{1}{\Gamma(1-\alpha)} \int_0^1 |g(\xi, v(\xi)) |\left[\int_0^T (T-\eta)^{-\alpha} C \eta^{\frac{\alpha}{2}-1} d\eta\right] d\xi + \frac{2C}{\alpha}T^{\frac{\alpha}{2}}\max_{[0,1]\times[0,T]}|f(x,t)|
$$
$$
\le CLT^{ - \frac{\alpha}{2}}  \frac{\Gamma\left(\frac{\alpha}{2}\right) }{\Gamma\left(1 - \frac{\alpha}{2}\right)} \int_0^1 | v(\xi)| d\xi + \frac{2C}{\alpha}T^{\frac{\alpha}{2}}\max_{[0,1]\times[0,T]}|f(x,t)|.
$$
Applying Grönwall's inequality, we obtain the following estimate:
$$
||v||\le  e^{\int_0^1LT^{ - \frac{\alpha}{2}} \cdot \frac{\Gamma\left(\frac{\alpha}{2}\right) }{\Gamma\left(1 - \frac{\alpha}{2}\right)}  d\xi }\frac{2C}{\alpha}T^{\frac{\alpha}{2}}\max_{[0,1]\times[0,T]}|f(x,t)|
= \frac{2C}{\alpha}T^{\frac{\alpha}{2}}e^{LT^{ - \frac{\alpha}{2}} \cdot \frac{\Gamma\left(\frac{\alpha}{2}\right) }{\Gamma\left(1 - \frac{\alpha}{2}\right)}}\max_{[0,1]\times[0,T]}|f(x,t)|.
$$
Therefore, if \( f(x,t) \in C([0,1] \times [0,T]) \), then \( F(x) \) and \( v(x) \) belong to \( C[0,1] \), which ensures that \( Av(x) \in C[0,1] \).

Next, we show that the operator \( A \) satisfies condition (2). To this end, again we successively make use of the properties of the operator $D^{\alpha-1}$, the upper bound for $|F(x)|$, the Lipschitz condition for the function $g$, and the estimate (\ref{NuraliyevaGrinasim}) for $G(x,T,\xi,0)$:

For any \( v_1(x), v_2(x) \in C[0,1] \), we estimate the following expression:

$$
|(Av_1)(x) - (Av_2)(x)|= \left| \int_0^1 D^{\alpha - 1}_{0T} G(x,T,\xi,0) [g(\xi, v_1(\xi)) - g(\xi, v_2(\xi))] \, d\xi \right|
$$

$$
=\left|\frac{1}{\Gamma(1-\alpha)} \int_0^1 \int_0^T (T-\eta)^{-\alpha} G(x,\eta,\xi,0) [g(\xi, v_1(\xi)) - g(\xi, v_2(\xi))] d\eta d\xi\right|
$$
$$
\le\frac{1}{\Gamma(1-\alpha)} \int_0^1 |g(\xi, v_1(\xi)) - g(\xi, v_2(\xi))|\left[\int_0^T (T-\eta)^{-\alpha} |G(x,\eta,\xi,0) | d\eta\right] d\xi
$$
$$
\le\frac{1}{\Gamma(1-\alpha)} \int_0^1 |g(\xi, v_1(\xi)) - g(\xi, v_2(\xi))|\left[\int_0^T (T-\eta)^{-\alpha}  C\eta^{\frac{\alpha}{2}-1} d\eta\right] d\xi
$$
$$
\le\frac{C}{\Gamma(1-\alpha)} \int_0^1 |g(\xi, v_1(\xi)) - g(\xi, v_2(\xi))|\left[T^{- \frac{\alpha}{2}}  \frac{\Gamma\left(\frac{\alpha}{2}\right) \Gamma(1 - \alpha)}{\Gamma\left(1 - \frac{\alpha}{2}\right)}\right] d\xi
$$
$$
\le CLT^{- \frac{\alpha}{2}} \frac{\Gamma\left(\frac{\alpha}{2}\right) }{\Gamma\left(1 - \frac{\alpha}{2}\right)} \int_0^1 | v_1(\xi)- v_2(\xi)| d\xi
\le \delta\max\limits_{x\in [0,1]} | v_1(x)- v_2(x)|
$$

$$
\delta=CLT^{ - \frac{\alpha}{2}} \frac{\Gamma\left(\frac{\alpha}{2}\right) }{\Gamma\left(1 - \frac{\alpha}{2}\right)}.
$$
According to the condition of Theorem~\ref{Nuraliyevatheo2}, we have \( \delta < 1 \).

Thus, \( A \) is a contraction operator in the space \( C[0,1] \). According to Banach's fixed point theorem, the operator \( A \) has a unique fixed point, i.e., the integral equation \eqref{Nuraliyevavx} has a unique solution \( v(x) \in C[0,1] \).

Given the non-local initial condition $u(x, 0) = g(x, u(x, T))$, and recognizing that $u(x, T) = v(x)$, we can now fully construct the solution $u(x, t)$ of the original problem. This is done using the representation \eqref{NuraliyevaGrin}, which expresses the solution via the Green function as:
\begin{equation}\label{Nuraliyevamainsol}
u(x, t) = \int_0^1  g(\xi, v(\xi)) D^{\alpha - 1}_{0t} G(x, t, \xi, 0)\, d\xi + \int_0^t \int_0^1 G(x, t, \xi, \tau) f(\xi, \tau)\, d\xi d\tau.    
\end{equation}

Therefore, under the assumptions stated in conditions \eqref{NuraliyevaCauchy1} and \eqref{NuraliyevaCauchy10}, we have constructed a well-defined, continuous, and explicit solution to the non-local and nonlinear boundary value problem involving a fractional Caputo derivative.

It is proven in \cite{NuraliyevaMamchuyev} that the function of the form~\eqref{Nuraliyevamainsol} satisfies the conditions of Definition~\ref{Nuraliyevadefen1}.


\end{document}